\documentclass[letterpaper, 10 pt, conference]{ieeeconf}  

\IEEEoverridecommandlockouts                              

\usepackage{graphicx}
\usepackage{algorithm}
\usepackage[noend]{algpseudocode}
\usepackage{amsmath, amssymb, mathtools}
\usepackage{caption}
\usepackage{subcaption}
\usepackage{booktabs}
\usepackage{comment}
\usepackage{tabularx}
\usepackage{xcolor}

\newtheorem{theorem}{Theorem}
\newtheorem{lemma}{Lemma}

\makeatletter
\newcommand{\multiline}[1]{%
  \begin{tabularx}{\dimexpr\linewidth-\ALG@thistlm}[t]{@{}X@{}}
    #1
  \end{tabularx}
}
\makeatother

\newcolumntype{Y}{>{\centering\arraybackslash}X}

\newcommand{\dom}{\prec}

\newcommand{\wdom}{\preceq}

\newcommand{\ndom}{\not \prec}

\newcommand{\alg}{RME-MOA*}

\title{\LARGE \bf
Runtime and Memory Efficient Multi-Objective A*
}

\author{Valmiki Kothare$^{1}$,
    Zhongqiang Ren$^{1}$, 
    Sivakumar Rathinam$^{2}$
    and
    Howie Choset$^{1}$%
\thanks{$^{1}$Valmiki Kothare, Zhongqiang Ren and Howie Choset are with the Robotics Institue, Carnegie Mellon University, 5000 Forbes Ave., Pittsburgh, PA 15213, USA {\tt\small  \{vok, zhongqir, choset\}@andrew.cmu.edu}}%
\thanks{$^{2}$Sivakumar Rathinam is with the Department of Mechanical Engineering, Texas A\&M University, College Station, TX 77843-3123, USA
        {\tt\small srathinam@tamu.edu}}%
}

\begin{document}

\maketitle
\thispagestyle{empty}
\pagestyle{empty}

\begin{abstract}
The Multi-Objective Shortest Path Problem (MO-SPP), typically posed on a graph, determines a set of paths from a start vertex to a destination vertex while optimizing multiple objectives. In general, there does not exist a single solution path that can simultaneously optimize all the objectives and the problem thus seeks to find a set of so-called Pareto-optimal solutions. To address this problem, several Multi-Objective A* (MOA*) algorithms were recently developed to quickly compute solutions with quality guarantees. However, these MOA* algorithms often suffer from high memory usage, especially when the branching factor (i.e. the number of neighbors of any vertex) of the graph is large. This work thus aims at reducing the high memory consumption of MOA* with little increase in the runtime. By generalizing and unifying several single- and multi-objective search algorithms, we develop the Runtime and Memory Efficient MOA* (RME-MOA*) approach, which can balance between runtime and memory efficiency by tuning two user-defined hyper-parameters.
\end{abstract}

\section{Introduction}\label{introduction}
Given a graph, the Shortest Path Problem (SPP) involves finding a solution path from a start vertex to a destination vertex with the least cumulative path cost.
As a natural extension of the SPP, the Multi-Objective Shortest Path Problem (MO-SPP) associates edges with cost vectors instead of scalar costs, where each element of the vector represents an objective to be optimized.
This problem arises in applications such as hazardous material transportation~\cite{erkut2007hazardous} where multiple objectives are optimized simultaneously. Instead, the MO-SPP attempts to find a set of Pareto-optimal (also called non-dominated) solution paths, and the corresponding set of cost vectors is often referred to as the Pareto-optimal front.
MO-SPP is computationally challenging as the number of Pareto-optimal paths can grow exponentially with respect to the size of the graph even with two objectives~\cite{hansen1980bicriterion}.

To compute the Pareto-optimal front, on the one hand, a family of MOA* algorithms were developed to quickly find the exact~\cite{NAMOA*dr,BOA*,BOBA*,EMOA*} or approximated~\cite{Zhang2022Apex,goldin2021approximate,warburton1987approximation} Pareto-optimal front, and this work limits its focus to the exact algorithms.
However, these algorithms primarily focus on the runtime efficiency and can consume a lot of memory during the search, especially in graphs with large branching factors (average number of successors per vertex), such as the state lattice graphs used in robotics~\cite{pivtoraiko2009differentially}, which limits their potential usage on platforms with limited computational resources (e.g. mobile robots, quadrotors).
On the other hand, while being able to provide theoretic guarantees on memory consumption, existing memory-efficient MOA* algorithms~\cite{PIDMOA*,IDMOA*} either make assumptions on the form of the input graph (e.g. acyclic graph), or sacrifice the runtime efficiency significantly, which hinders their practical usage in large graphs or dense state lattices in robotics~\cite{pivtoraiko2009differentially}.
This work thus develops a new algorithm, called Runtime- and Memory-Efficient MOA* (\alg) that balances runtime and memory efficiency by leveraging several techniques from previous works.

\alg{} has three distinct features:
Firstly, its search framework is based on our prior EMOA*~\cite{EMOA*}, a runtime-efficient MOA* algorithm
Secondly, it extends the \emph{partial expansion} technique~\cite{PEA*} from the single- to multi-objective case and fuses it with EMOA* in order to trade off runtime for memory efficiency.
Finally, \alg{} leverages PIDMOA*~\cite{PIDMOA*}, a memory-efficient iterative deepening depth-first MOA* algorithm, by switching between best-first search, as in EMOA*, and iterative depth-first search, as in PIDMOA*, to further improve memory efficiency.
We show that the proposed \alg{} is a strict generalization of EMOA* and PIDMOA* and can move along a spectrum from runtime efficient (EMOA*) to memory efficient (PIDMOA*) by tuning two user-defined hyper-parameters that control the search process.

Our \alg{} is guaranteed to find the Pareto-optimal solution front.
We test our \alg{} in both grid-like graphs with varying branching factors and robot state lattices against two baselines: its predecessors EMOA* and PIDMOA*.
\alg{} is able to significantly reduce memory consumption compared to EMOA*, at the expense runtime performance, on problem instances with large branching factors.
\alg{} is also able to significantly increase runtime performance compared to PIDMOA* at the expense of memory efficiency.

\section{Problem Description}\label{problem}
Let $G = (V,E,\vec{c})$ denote a directed graph, where $V$ is a set of vertices, $E$ a set of directed edges connecting those vertices, and $\vec{c}: E \to (\mathbb{R}^{\geq 0})^M$ a cost function that maps edges to non-negative cost vectors, where $M$ is the number of objectives to be minimized. Let $Succs(v) = \{v'\in V| \exists (v,v') \in E\}$ be the successor vertices in graph $G$ of vertex $v$.

Let $\pi(v_1, v_n)$ be a path in $G$ from $v_1$ to $v_n$, which is a sequence of vertices $\{v_1, v_2, ... , v_n\}$ such that $\exists (v_i, v_{i+1})\in E$ for all $i \in \{1,...,n-1\}$.
Let $v_s$ and $v_d$ denote the start and destination vertex respectively.
We call a path $\pi(v_s, v)$ from the start vertex $v_s$ to some arbitrary vertex $v \neq v_d$ a \textit{partial path}.
Every path $\pi(v_1, v_n)$ is associated with a vector path cost $\vec{g}(\pi(v_1, v_n)) = \sum_{i=1}^{n-1} \vec{c}(v_i, v_{i+1})$ where $(v_i, v_{i+1})\in E$.
To compare two paths, we compare their corresponding cost vectors using the notion of dominance \cite{MOA*}.
Given $\vec{a},\vec{b} \in (\mathbb{R}^{\geq 0})^M$ with $a_i$ denoting the $i$-th component of $\vec{a}$, $\vec{a}$ \textit{dominates} $\vec{b}$ ($\vec{a} \dom \vec{b}$)%
\footnote{In the literature, $a$ dominates $b$ and $a$ weakly dominates $b$ are sometimes also denoted as $a \succeq b$ and $a \leq b$, respectively \cite{EMOA*}. We choose to use notation $a \dom b$ and $a \wdom b$ (as in \cite{PIDMOA*}) as it reads more intuitively and alleviates confusion with lexicographic ordering ($\leq_{lex}$).}
iff $\forall i,\: a_i \leq b_i \; \land \; \vec{a} \neq \vec{b}$.
Similarly, $\vec{a}$ \textit{weakly dominates} $\vec{b}$ ($\vec{a} \wdom \vec{b}$) iff $\forall i, a_i \leq b_i$.

Let $\pi(v_1, v_n)$ and $\pi'(v_1, v_n)$ be two distinct paths from $v_1$ to $v_n$. $\pi$ dominates $\pi'$ ($\pi \dom \pi'$) iff $\vec{g}(\pi) \dom \vec{g}(\pi')$, and $\pi$ weakly dominates $\pi'$ ($\pi \wdom \pi'$) iff $\vec{g}(\pi) \wdom \vec{g}(\pi')$.
Note in the case where $v_1 = v_s$ and $v_n = v_d$, the paths are solution paths from the start vertex to the destination vertex.
In the MO-SPP, a \textit{Pareto-optimal solution} is a path $\pi^*(v_s,v_d)$ that is not dominated by any other path $\pi'(v_s,v_d)$ in $G$ (i.e. $\pi' \ndom \pi^*$ for all $\pi' \neq \pi^*$ in $G$).
The set of all such paths $\pi^*$ is the exact Pareto-optimal (or non-dominated) solution set.
Any maximal subset $\mathcal{S}$ of the Pareto-optimal solution set where any two paths in $\mathcal{S}$ do not have the same cost is known as a maximal \textit{cost-unique Pareto-optimal} solution set.
The goal of this work is to find a maximal cost-unique Pareto-optimal solution set for the MO-SPP.

\section{Preliminaries}\label{prelim}
\subsection{Notation and Terminology}\label{prelim:notation}
Let $\vec{h}: V \to (\mathbb{R}^{\geq 0})^M$ be a heuristic function mapping vertices to heuristic vectors, where $\vec{h}(v)$ is an under-estimate of the true optimal path cost from $v$ to $v_d$.
Our algorithm, as well as its predecessors \cite{EMOA*, BOA*, NAMOA*dr} rely on the heuristic being consistent:
$\vec{h}$ is consistent iff
$\vec{h}(v_d) = \vec{0} \; \land \; \forall (v,v')\in E,\: \vec{h}(v) \leq \vec{c}(v, v') + \vec{h}(v')$.
This guarantees that $\vec{h}$ does not overestimate the cost between any two adjacent vertices.
Let $\vec{f}(\pi(v_s,v)) := \vec{g}(\pi(v_s,v)) + \vec{h}(v)$ denote the lower bound on the vector cost of any solution path that proceeds from partial path $\pi(v_s,v)$.

Let a \textit{label} $l$ represent a partial path $\pi(v_s,v)$.%
\footnote{To identify a partial path, different names such as nodes~\cite{BOA*}, states~\cite{ren22mopbd} and labels~\cite{martins1984multicriteria,sanders2013parallel} have been used in the multi-objective path-planning literature. This work uses ``labels'' to identify partial solution paths.}
A label $l$ is \emph{expanded} when the corresponding partial path is extended to all successors of $v(l)$.
A label is dominated (or weakly dominated) iff the corresponding partial path $\pi(l)$ is dominated (or weakly dominated).
Let $v(l) = v$, $\vec{g}(l) = \vec{g}(\pi)$, $\vec{f}(l) = \vec{f}(\pi)$, and $parent(l)$ be the parent label that was expanded to create $l$.
Let $Child(l) = \{l' | v(l') = v', \forall v' \in Succs(v(l))$ be the children labels of parent label $l$ that are created on expansion of $l$.
Let OPEN be an ordered list of labels, implemented as a priority queue, where labels are ordered lexicographically by $\vec{f}$.
Let CLOSED denote a set of labels that have been expanded.
Let $\alpha(v)$ be the \textit{Pareto frontier} at vertex $v$, i.e. the set of all non-dominated partial path costs $\vec{g}$ (and their associated labels) found so far during the MOA* search from $v_s$ to $v$.
We say $\alpha(v)$ dominates (or weakly dominates) $\vec{g}(l)$, iff $\exists l' \in \alpha(v)$ s.t. $\vec{g}(l') \dom \vec{g}(l)$ (or $\vec{g}(l') \wdom \vec{g}(l)$), which is denoted as $\alpha(v) \dom \vec{g}(l)$ (or $\alpha(v) \wdom \vec{g}(l)$)).
Let $\mathcal{S}$ be the Pareto-optimal solution frontier, which, when the search terminates, is a maximal cost-unique Pareto-optimal solution set.

\subsection{EMOA*}\label{prelim:EMOA*}
Because we base \alg{} on EMOA*'s search framework, we elaborate EMOA* here, referring to Algorithms 1,2, and 3 in ~\cite{EMOA*} to save space.
At first, a label for the start vertex is created and added to OPEN.
The main loop then iteratively \textit{extracts} from OPEN a label $l$ with the lexicographically smallest $\vec{f}(l)$.
The extracted label $l$ is then \emph{checked} for dominance: if $\alpha(v(l)) \wdom \vec{g}(l)$ or $\mathcal{S} \wdom \vec{f}(l)$, which are called \emph{FrontierCheck} and \emph{SolutionCheck}, respectively (Alg.~2 in~\cite{EMOA*}). 
If $l$ is not weakly dominated by either $\alpha(v(l))$ or $\mathcal{S}$, then, a procedure called \emph{UpdateFrontier} is invoked (Alg.~3 in~\cite{EMOA*}), where $l$ is added to $\alpha(v(l))$, and all labels in $\alpha(v(l))$ whose $\vec{g}$ are dominated by $\vec{g}(l)$ are \emph{filtered}.
Then, $l$ is expanded, and for each newly generated label $l'$ dominance checks are performed to determine if either $\alpha(v') \wdom \vec{g}(l')$ or $\mathcal{S} \wdom \vec{f}(l')$.
If not, then $l'$ is added to OPEN.
The search terminates when OPEN is empty, which guarantees that $\mathcal{S} = \alpha(v_d)$ is a maximal cost-unique Pareto optimal solution set.

The frequent dominance checks during MOA* search pose a computational challenge, which EMOA* addresses in three ways. 
Firstly, it uses ``dimensionality reduction''~\cite{NAMOA*dr} to reduce the number and dimension of cost vector comparisons in each \emph{Check} and \emph{Filter} operation~\cite{EMOA*}. These operations are additionally sped up by using a balanced binary search tree (BBST) to represent $\alpha(v)$ at each $v\in V$, which further reduces the number of vector comparisons per operation. Finally, EMOA* leverages the notion of ``lazy checks'' from~\cite{BOA*} to defer the dominance checks of a label until just before its expansion, which removes the need to explicitly filter OPEN after each label generation. We refer the reader to~\cite{EMOA*} for a more detailed explanation of these techniques. While these improvements help in improving runtime performance, EMOA* still suffers from high memory consumption.

\subsection{PIDMOA*}\label{prelim:PIDMOA*}
Unlike EMOA*, which addresses the high runtime expense of the MO-SPP, PIDMOA*~\cite{PIDMOA*} attempts to do so for memory consumption. 
An extension of Iterative Deepening A* (IDA*)~\cite{IDA*} to the multi-objective case, PIDMOA* uses the same strategy of performing DFS up to a threshold and iteratively increasing the threshold after each search iteration until the optimal solution is found. 
Instead of a scalar cost threshold as in IDA*, however, PIDMOA* maintains a Pareto threshold set ($\mathcal{T}$) of cost vectors. 
Search is terminated at any partial path $\pi$ where $\mathcal{T} \dom \vec{f}(\pi)$ or $\mathcal{S} \wdom \vec{f}(\pi)$, where, just as in EMOA*, $\mathcal{S}$ is the Pareto-frontier maintained at the destination vertex. 
PIDMOA* is designed to consume memory linear to the depth and number of the solutions, as DFS only explicitly maintains the current path (and all found solutions) during search.
Notably, PIDMOA* does not maintain a CLOSED set, as EMOA* does in the form of $\alpha$.
As a result, however, its runtime increases dramatically with the depth of the solutions. 
Refer to~\cite{PIDMOA*} for a more detailed explanation of the original algorithm.

\subsection{PEA*}\label{prelim:PEA*}
Our method leverages the partial expansion technique described in \cite{PEA*}.
Different from A*, Partial Expansion A* (PEA*) only explores a \emph{subset} of successors of $v$ when expanding $v$ and reserves the remaining successors for future expansion by re-inserting their predecessor $v$ into OPEN with an ``updated $f$ value'', which we call its \emph{re-expansion} value $r(v)$.%
\footnote{In \cite{PEA*}, $r(v)$ of a vertex $v$ is called the \textit{stored} value of $n$ and denoted $F(v)$. We change notation to $r(v)$ and the name to \textit{re-expansion} value to alleviate potential confusion between $f(v)$ and $F(v)$ and clarify the purpose of the value.}
To expand a vertex $v$, PEA* explores only the subset of the successors $v'$ of $v$ that satisfy $f(v') \leq r(v) + C$ with $C \in \mathbb{R}^{\geq 0}$ being a scalar constant hyperparameter that parameterizes the ``amount'' of partial expansion, which we explain later.
All successors that satisfy this condition are \textit{explored} and are expanded as in regular A*.
The other successors that do not satisfy this condition are \textit{unexplored} and are reserved for future \textit{re-expansion} by re-inserting $v$ into OPEN (now ordered by $r$ instead of $f$) with an updated $r(v) = \min f(v')$ for all unexplored successors $v'$ of $v$.
$v$ is later re-expanded one or more times until all its successors have been explored.
PEA* is thus able to trade off runtime (increased number of (re)expansions) for memory efficiency (less explored successors per (re)expansion, resulting in a smaller maximal size of OPEN).
This tradeoff is parameterized by $C$, where small values of $C$ correspond to less explored successors per (re)expansion and more (re)expansions overall, and sufficiently large values behave identically to A*.
PEA* has been shown to be memory efficient for large branching factors problems, claiming to reduce memory consumption by a factor of the branching factor $b$.
This is because PEA* reduces the size of OPEN by, at most, a factor of $b$ because the ratio of OPEN to CLOSED is, at most, $b-1$ to $1$ in traditional A*~\cite{PEA*}.
Notably, PEA* does not reduce the size of the CLOSED set and so is less effective on small branching factor problems, where the ratio of OPEN to CLOSED is small.

\section{Method}\label{method}
\begin{algorithm}[tb]
\small
\caption{\alg}\label{alg:\alg}
\begin{algorithmic}[1]
\State \multiline{%
$l_s \gets$ new label: $v(l_s) \gets v_s$, $\vec{g}(l_s) \gets \vec{0}$, $\vec{f}(l_s) \gets \vec{h}(v_s)$, $\vec{r}(l_s) \gets \vec{h}(v_s)$}
\State Add $l_s$ to OPEN
\State $\alpha (v) \gets \emptyset, \forall v \in V$
\State $\vec{C} \in (\mathbb{R}^{\geq 0})^M$ \Comment{User-Defined}
\While {OPEN $\neq \emptyset$}
    \State $l \gets$ OPEN.pop() \Comment{Ordered lexicographically}
    \If {\emph{FrontierCheck($l$)} \textbf{or} \emph{SolutionCheck($l$)}}
        \State \textbf{continue}
    \EndIf
    \If {$v(l) = v_d$}
        \State \emph{UpdateSolution($\mathcal{S}$, $l$)}
        \State \textbf{continue} 
    \EndIf
    \State \emph{UpdateFrontier($l$)}
    \If {$\vec{h}(v(l)) < \vec{D}$}
        \State \emph{PIDMOA*}($v(l)$, $v_d$, $\vec{g}(l)$)
        \State \textbf{continue}
    \EndIf
    \State $\vec{r}_{next}(l) \gets \vec{\infty}$
    \For {\textbf{all} $v' \in Succs(v(l))$} \Comment{Label (Re)expanded}
        \State \multiline{%
        $l' \gets$ new label: $v(l') \gets v'$, $\vec{g}(l') \gets \vec{g}(l) + \vec{c}(v, v')$, $\vec{f}(l') \gets \vec{g}(l') + \vec{h}(v(l'))$, $\vec{r}(l') \gets \vec{g}(l') + \vec{h}(v(l'))$}
        \State $parent(l') \gets l$
        \If {$\vec{f}(l') <_{lex} \vec{r}(l)$}
            \State \textbf{continue}
        \EndIf
        \If {\emph{FrontierCheck($l'$)} \textbf{or} \emph{SolutionCheck($l'$)}}
            \State \textbf{continue}
        \EndIf
        \If {$\vec{f}(l') >_{lex} \vec{r}(l) + \vec{C}$} \Comment{Unexplored Child}
            \State $\vec{r}_{next}(l) \gets \min_{lex}(\vec{r}_{next}(l), \vec{f}(l'))$
            \State \textbf{continue}
        \EndIf
        \State Add $l'$ to OPEN \Comment{Explored Child}
    \EndFor
    \If {$\vec{r}_{next}(l) \neq \vec{\infty}$}
        \State $\vec{r}(l) \gets \vec{r}_{next}(l)$
        \State Add $l$ to OPEN
    \EndIf
\EndWhile
\State return $\mathcal{S}$
\end{algorithmic}
\end{algorithm}

\begin{algorithm}[tb]
\small
\caption{\emph{PIDMOA*}($v_p$, $v_d$, $g_p$)}\label{alg:PIDMOA*}
\begin{algorithmic}[1]
\State $l_p \gets$ new label: $v(l_p) \gets v_p$, $\vec{g}(l_p) \gets g_p$, $\vec{f}(l_p) \gets g_p + \vec{h}(v_p)$
\State $\mathcal{T} \gets \{\vec{h}(v(l_p))\}$
\While {$\mathcal{T} \neq \emptyset$}
    \State Add $l_p$ to OPEN\_DFS
    \State $\mathcal{T}_n = \emptyset$
    \While {OPEN\_DFS $\neq \emptyset$}
        \State $l \gets$ OPEN\_DFS.pop()
        \If {\emph{SolutionCheck($l$)}}
            \State \textbf{continue}
        \EndIf
        \For {\textbf{all} $v' \in Succs(v(l))$}
            \State \multiline{%
            $l' \gets$ new label: $v(l') \gets v'$, $\vec{g}(l') \gets \vec{g}(l) + \vec{c}(v, v')$, $\vec{f}(l') \gets \vec{g}(l') + \vec{h}(v(l'))$}
            \State $parent(l') \gets l$
            \If {\emph{SolutionCheck($l'$)}}
                \State \textbf{continue}
            \EndIf
            \If {\emph{ThresholdCheck($l'$)}}
                \State \emph{UpdateThreshold($\mathcal{T}_n$, $l'$)}
                \State \textbf{continue}
            \EndIf
            \If {$v(l') = v_d$}
                \State \emph{UpdateSolution($\mathcal{S}$, $l'$)}
                \State \textbf{continue}
            \EndIf
            \State Add $l'$ to OPEN\_DFS
        \EndFor
    \EndWhile
    \State $\mathcal{T} \gets \mathcal{T}_n$
\EndWhile
\end{algorithmic}
\end{algorithm}


\subsection{New Concepts}\label{method:notation}
Let $\vec{r}(l)$ be the $M$-dimensional \emph{re-expansion} vector associated with a label $l$.
As with $r(v)$ in PEA*, $\vec{r}(l)$ is used to efficiently ``store'' information about the children of $l$ to determine when it is necessary to re-expand $l$.
Let $<_{lex}, and >_{lex}$ be lexicographic comparison operators between two vectors.
On label generation, $\vec{r}(l) = \vec{f}(l)$, and at the end of each (re)expansion, $\vec{r}(l) = \min_{lex}\vec{f}(l')$ for all \emph{unexplored} children $l'$, defined as follows: 
For each (re)expansion of $l$, let $\mathcal{E}(l) = \{l' \in Child(l)|\vec{r}(l) <_{lex} \vec{f}(l') \leq_{lex} \vec{r}(l) + \vec{C}\}$ be the \emph{explored} children of $l$.
Let $\mathcal{U}(l) = \{l' \in Child(l)|\vec{f}(l') >_{lex} \vec{r}(l) + \vec{C}\}$ be the \textit{unexplored} children of $l$.%
\footnote{In \cite{PEA*}, these explored and unexplored sets of children are called \textit{promising} and \textit{unpromising}, respectively. We change semantics here because the relative `` promise'' of one label vs another is not obvious for multi-objective costs and lexicographic ordering does not fully capture this notion.}
Let $\vec{C} \in (\mathbb{R}^{\geq 0})^M$ be a user-defined constant vector that parameterizes the degree of partial expansion, as in PEA*.
Let $\mathcal{P}(l) = \{l' \in Child(l)|\vec{f}(l') \leq_{lex} \vec{r}(l)\}$ be the \emph{previously explored} children of $l$, which will always be $\emptyset$ on first expansion of $l$.

Let $\vec{D} \in (\mathbb{R}^{\geq 0})^M$ be another user-defined constant vector that parameterizes the ``distance'' from the goal vertex at which to start PIDMOA* search, where distance from $v$ to $v_d$ is estimated by $\vec{h}(v)$
(to avoid confusion, we call the start vertex for each PIDMOA* iteration $v_p$ instead of $v_s$).
Let OPEN\_DFS be a stack that stores labels in Last-In First-Out (LIFO) order. 
Let $\mathcal{T}$ be the current Pareto threshold frontier maintained during each PIDMOA* search iteration, and let $\mathcal{T}_n$ be the Pareto threshold frontier for the next search iteration.

For all $l \in \mathcal{S}$, let $\mathcal{S}_{path}(l)$ be the \emph{PIDMOA* path} of $l$, or the portion of the path $\pi(l)$ found during PIDMOA* search.

\subsection{Algorithm Overview}\label{method:alg_overview}
Alg.~\ref{alg:\alg} depicts the pseudocode for \alg, and it follows a similar framework to that described in Sec.~\ref{prelim:EMOA*}.
In each \alg{} iteration (lines 5-31), label $l$ with the lexicographically smallest $\vec{r}$ is extracted from OPEN (line 6). As in EMOA*, \emph{FrontierCheck} and \emph{SolutionCheck} are performed, but now $v(l) = v_d$ is checked before deciding whether to call the new operation \emph{UpdateSolution} or \emph{UpdateFrontier} (lines 7-12).
\emph{FrontierCheck} and \emph{UpdateFrontier} are additionally modified such that, if $l$ has been added to $\alpha(v(l))$ before, \emph{FrontierCheck} returns false, and \emph{UpdateFrontier} does not modify $\alpha(v(l))$
These modifications are discussed in Sec.~\ref{method:discussion:frontier_ops}. 
On lines 13-14, if $v(l)$ is within a certain distance of $v_d$ for all cost functions (i.e. $\vec{h}(v(l)) < \vec{D}$), $l$ is now expanded using a modified PIDMOA* search, depicted in Alg.~\ref{alg:PIDMOA*} (discussion in Sec.~\ref{method:discussion:pidmoa_discussion}). 
This is done in an effort to save memory while avoiding the steep runtime cost of a deep PIDMOA* search.
Otherwise, $l$ is expanded using partial expansion. 
For all $l' \in Child(l)$, only those that satisfy what we call the \textit{partial expansion check} (lines 18-27) are explored and added to the OPEN (discussed in Sec.~\ref{method:discussion:pe_check}).
If, after expansion, $l$ has any remaining unexplored children, $l$ is added back to OPEN with $\vec{r}(l) = \min_{lex}\vec{f(l')}$ for re-expansion (lines 16, 25, 28-30).
The search process iterates until OPEN depletes.
At termination, $\mathcal{S}$ is returned (line 31), which is guaranteed to be a maximal cost-unique Pareto-optimal solution set (Sec.~\ref{theoretical_analysis}).

\subsection{Key Procedures and Discussion}\label{method:discussion}
\subsubsection{Partial Expansion Check}\label{method:discussion:pe_check}

For each expansion of $l$, the partial expansion check (lines 18-27) determines if $l' \in Child(l)$ should be added to OPEN or ``saved'' for later re-expansion of $l$.
As in EMOA*, $l'$ must pass \emph{FrontierCheck($l'$)} and \emph{SolutionCheck($l'$)} in order to be considered for exploration (line 22). 
All $l' \in \mathcal{U}(l)$ are used to update $\vec{r}_{next}(l)$ and discarded, while the remaining $l' \in \mathcal{E}(l)$ are added to OPEN.
We choose to use lexicographic comparison on lines 20, 24, and 25 in order to maintain the optimality guarantees of EMOA*: by setting $\vec{r}(l) = \min_{lex}\vec{f(l')}$ for all $l' \in \mathcal{U}(l)$ before adding $l$ back to OPEN, the algorithm guarantees that all unexplored children will eventually be re-expanded when or before they are needed. 
In subsequent re-expansions of $l$, all $l' \in \mathcal{P}(l)$ are immediately discarded (line 20) as they have already been explored.
This, too, is guaranteed by the use of lexicographic comparison because $\forall l' \in \mathcal{E}_{n-1}(l), \vec{f}(l') < \vec{r}_{n-1}(l) + \vec{C} < \vec{r}_{n}(l)$, where $n$ refers to the current number of (re)expansions of $l$. 
Thus $\mathcal{P} = \bigcup_{i=1}^{n-1} \mathcal{E}_i$, so it is safe to discard any $l' \in \mathcal{P}(l)$ because it will have already been added to OPEN or have failed \emph{FrontierCheck} or \emph{SolutionCheck} and have been previously discarded.
Distinct from PEA*, the dominance checks (line 22) are performed before the ``exploration'' check (line 24). 
This modification prevents unnecessary re-expansions of $l$ and results in better computational and memory performance in experimental testing.

\subsubsection{Modified Frontier Operations}\label{method:discussion:frontier_ops}
For \alg{} to be guaranteed to produce a maximal cost-unique Pareto-optimal solution set, we must replace $\alpha(v_d)$ with a separate $\mathcal{S}$ that does not utilize the dimensionality reduction proposed in~\cite{EMOA*}. 
Unlike EMOA*, PIDMOA* does not guarantee that the sequence of labels being expanded at $v_d$ has non-decreasing $f_1$ values.
Because some solution paths may be reached by partial expansion and others by PIDMOA* search, we now cannot use this non-decreasing $f_1$ assumption to reduce the dimension of $\alpha(v_d)$, and thus must use the full label vectors for \emph{SolutionCheck} and \emph{UpdateSolution}, which operate on this new $\mathcal{S}$ and not $\alpha(v_d)$.
$\mathcal{S}_{path}$ additionally stores the PIDMOA* for each solution, which is necessary to rebuild the path back from $v_d$ to $v_p$ (the start vertex for PIDMOA* search), after which $\alpha$ can be used to rebuild the path from $v_p$ to $v_s$.

We also modify \emph{FrontierCheck} and \emph{UpdateFrontier} from EMOA* as follows: if $l$ has been previously added to $\alpha(v(l))$, \emph{FrontierCheck($l$)} returns false and \emph{UpdateFrontier($l$)} returns without modifying $\alpha(v(l))$. 
Otherwise, \emph{FrontierCheck} and \emph{UpdateFrontier} proceed as in EMOA*.
These are necessary modifications due to the introduction of partial expansion, which will cause \emph{FrontierCheck} and \emph{UpdateFrontier} to be called on $l$ multiple times and would otherwise stop the re-expansion of $l$ without these changes (on re-expansion, $l \in \alpha(v(l))$, therefore $\alpha(v(l)) \wdom l$, so \emph{FrontierCheck($l$)} would return false).

We finally add two procedures \emph{ThresholdCheck} and \emph{UpdateThreshold}, which perform the dominance check $\mathcal{T} \dom \vec{f}(l)$ and the dominance filter and update of $\mathcal{T}_n$ using $\vec{f}(l)$ for some label $l$, respectively.
These procedures are nearly identical to \emph{SolutionCheck} and \emph{UpdateSolution} (implemented using BBSTs without dimensionality reduction) except that \emph{ThresholdCheck} uses dominance instead of weak dominance.


\subsubsection{Modified PIDMOA* Search}\label{method:discussion:pidmoa_discussion}
We modify the approach in~\cite{PIDMOA*} by implementing PIDMOA* iteratively (depicted in Alg.~\ref{alg:PIDMOA*}) instead of recursively due to the runtime detriment of recursive implementation.
At each extraction of $l$ from OPEN, a PIDMOA* search instance is started if $\vec{h}(v(l)) < \vec{D}$.
A new label $l_p$ with $v(l_p) = v_p = v(l)$ and $\vec{g}(l_p) = g_p = \vec{g}(l)$ is added to OPEN\_DFS at the start of each search iteration. 
$\mathcal{T}$ is initialized to $\{\vec{h}(v_p)\}$ and at the beginning of each search iteration $\mathcal{T}_n$ is initialized to $\emptyset$.
Labels are iteratively extracted from OPEN\_DFS in LIFO order, and if a label $l$ passes \emph{SolutionCheck} (line 8), it is expanded.
If a child $l' \in Child(l)$ does not satisfy \emph{SolutionCheck}, it is discarded. 
$l'$ is similarly discarded if it does not satisfy \emph{SolutionCheck($l'$)} but is first used in \emph{UpdateThreshold($\mathcal{T}_n$, $l'$)}.
This iterative construction of $\mathcal{T}_n$ from the partial paths immediately dominated by $\mathcal{T}$ (followed by the discontinuation of the current search iteration at those paths) guarantees the optimality of PIDMOA* and allows for the representation of the threshold sets as BBSTs to improve runtime for dominance checks and filtering.
All non-discarded children $l'$ are then added to OPEN\_DFS.
A single search iteration ends when OPEN\_DFS is empty, after which $\mathcal{T}$ is set to $\mathcal{T}_n$ and a new iteration begins.
This continues until $\mathcal{T}$ is empty after an iteration, guaranteeing that a maximal cost-unique Pareto-optimal solution set extending from $l_p$ is found, conditioned on the previously found non-dominated solutions and $g_p$.
The use of BBSTs for $\mathcal{T}$ and $\mathcal{S}$ in this modified PIDMOA*, coupled with the iterative implementation, significantly reduced the runtime of PIDMOA* in our tests when compared to the original approach~\cite{PIDMOA*}.

\subsubsection{Relationship to Predecessors}\label{method:discussion:predecessors}
The intuition behind the combined use of EMOA* with partial expansion and PIDMOA stems from the fact that partial expansion reduces the size of OPEN while PIDMOA* reduces (eliminates) the size of CLOSED.
With a $\vec{D} = \vec{0}$, \alg{} strictly generalizes EMOA* with the partial expansion technique: that is, by setting $\vec{C} = \vec{\infty}$ in \alg, we implement the same behavior as EMOA*. With a non-zero $\vec{D}$, \alg{} also strictly generalizes PIDMOA* (with our implementation-specific modifications). This is done by setting $\vec{D} = \vec{\infty}$, causing \alg{} to behave identically to PIDMOA*.

\section{Theoretical Analysis of \alg}\label{theoretical_analysis}
\begin{figure*}[t]
    \centering
    \begin{tabular}{c|c}
    \toprule
    \multicolumn{2}{c}{20x20, $2^k$-connected, M=2 Grid Test}\\
    \midrule
    $k = 2$ & $k = 5$\\
    \includegraphics[height=0.22\textwidth]{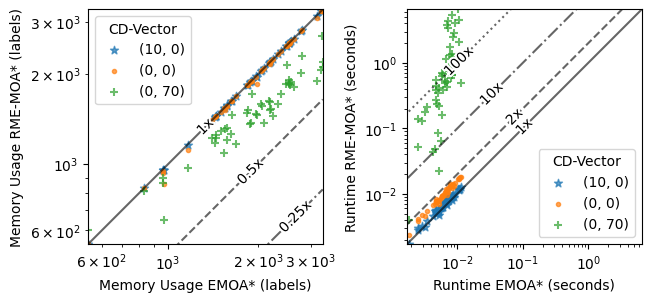} & \includegraphics[height=0.22\textwidth]{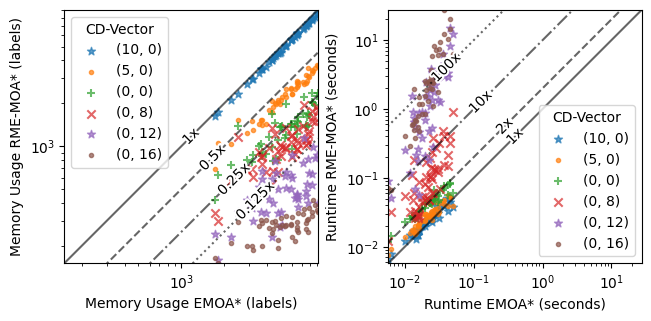}\\
    \end{tabular}
    \caption{Comparison of memory usage and runtime metrics for \alg{} with varying tuples $(\vec{C}, \vec{D})$, against EMOA* $(\vec{C}=\vec{\infty}, \vec{D}=\vec{0})$. The search space is a 20x20 empty grid environment of $2^k$-connectedness and $M=2$. For each of 50 randomized grid instances, we run \alg{} with all $(\vec{C}, \vec{D})$ tuples and EMOA*. The lines on graphs represent multipliers of the performance of EMOA* for each metric (e.g. in (c), any point below the ``0.5x'' line represents a search instance of \alg{} that had less than half the memory usage of the corresponding search instance of EMOA*). Plots are in log scale. Note that search instances of EMOA* are not expressly plotted as they would all lie along the ``1x'' line. Here, we compare performance of \alg{} on varying branching factors, with $k=2$ on the left and $k=5$ on the right. Note that the specific $(\vec{C}, \vec{D})$ tuples differ between branching factors: this is because the choice of $\vec{C}$ and $\vec{D}$ depend on the branching factor and maximum depth of a solution, both of which vary between different $2^k$-connected grids.}
    \label{fig:grid}
\end{figure*}
(Note: for Lemmas 1 and 2 and Theorem 1, we first consider \alg{} with $\vec{D} = 0$ (i.e. no PIDMOA* search). For the remaining Lemma and Theorem, we  treat $\vec{D}$ as non-zero)
\begin{lemma}\label{lemma:1}
    At any time during the search, let $l_{min}$ denote the lexicographically smallest label in OPEN. For any label $l \in$ OPEN, all of its unexplored children $l' \in \mathcal{U}(l)$ must satisfy $\vec{f}l_{min} \leq_{lex} \vec{f}(l')$.
\end{lemma}
\begin{lemma}\label{lemma:2}
    \alg{} (with $\vec{D} = 0$) and EMOA* expand labels in the same order.
\end{lemma}
\begin{theorem}\label{theorem:1}
    \alg{} (with $\vec{D} = 0$) computes a maximal cost-unique Pareto-optimal solution set.
\end{theorem}
\begin{lemma}\label{lemma:3}
    Given a label $l_p$ representing some partial path to vertex $v_p$ with cost $g_p$, PIDMOA* search returns a maximal cost-unique Pareto-optimal solution set extending $l_p$ from $v_p$ to $v_d$.
\end{lemma}
\begin{lemma}\label{lemma:4}
    All paths $\pi$ in the exact Pareto-optimal solution set are reached by \alg{}
\end{lemma}
\begin{theorem}\label{theorem:2}
    \alg{} (for any $\vec{D} \in (\mathbb{R}^{\geq 0})^M$) computes a maximal cost-unique Pareto-optimal solution set.
\end{theorem}

\section{Numerical Results}\label{evaluation}
\subsection{Baselines and Metrics}\label{evaluation:metrics}
Because \alg{} strictly extends both EMOA* and PIDMOA* (see Sec.~\ref{method:discussion:predecessors}, we compare \alg{} to EMOA* and plot PIDMOA* results (where computationally feasible).
In our tests, we compare runtime and memory usage between a spectrum of $\vec{C}$ and $\vec{D}$ values.
Runtime is measured as the time taken (in seconds) to find a maximal cost-unique Pareto-optimal solution, not including time to construct the graph and calculate the heuristic, which is negligible.
Memory usage is quantified by the maximum number of \textit{stored} labels, defined as:
\begin{equation*}
\begin{split}
     &\max_t(|\text{OPEN}_t| + |\text{OPEN\_DFS}_t| + |\mathcal{T}_t| + \\
     &|(\mathcal{T}_n)_t| + \sum_{v\in V}|\alpha_t(v)| + \sum_{l \in \mathcal{S}}|(\mathcal{S}_{path})_t(l)|)
\end{split}
\end{equation*}
for all $t \in [0,t_e]$, where $t$ is some time during the search execution, $t_e$ is the time at the end of execution, and the $|S_t|$ operation represents the number of labels stored in some set $S$ at time $t$.%
\footnote{For simplicity, for a label $l \in \mathcal{S}$ reached by partial expansion, we say $|(\mathcal{S}_{path})_t|(l) = 1$ and so forego the need to count $|\mathcal{S}_t|$}
This captures the number of labels necessary to find the Pareto-optimal solution set and to reconstruct the solution paths back from the destination vertex at the end of execution.%
\footnote{We inherit this use of stored labels for measuring memory usage from prior literature~\cite{PEA*, EPEA*}. Being an implementation- and hardware-independent metric, it offers insight into the reason behind changes in memory usage as result of parameter changes. We recognize that some memory overhead is incurred independent of label storage (i.e. heuristic storage), and also that it may not be necessary to store all elements of a label in OPEN\_DFS, $\alpha$, $\mathcal{T}$, $\mathcal{T}_n$, or $\mathcal{S}_{path}$.}

We implement all algorithms in C++ and test on a Windows 11 laptop running WSL2 Ubuntu 20.04 with an AMD Ryzen 9 6900HS 3.30GHz
CPU and 16 GB RAM without multi-threading. As in~\cite{EMOA*}, for each component $c_i$ of the cost objective function $\vec{c}$, the heuristic is pre-computed using Djikstra's search from $v_d$ to all $v \in V$ with edge costs $c_i(e), \forall e \in E$.

\subsection{Grid Map Test}\label{evaluation:grid_map}
We first compare \alg{} with varying $\vec{C}$ and $\vec{D}$ in an empty \textit{$2^k$-connected grid}, where each cell in the grid has $2^k$ neighbors (see~\cite{MAPF} for a visualization).
This allows us to parameterize the branching factor of the graph via $k$ and thus compare the performance of \alg{} on graphs with different branching factors. 
We formulate a fully parameterized testing framework where $M$ is the number of objectives, $\vec{C}$ the partial expansion parameter, $\vec{D}$ the PIDMOA* search parameter, $k$ the degree of connectedness, and $r \times c$ the dimensions of the grid. 
Each objective in the vector cost function is uniformly randomly sampled from $\{1,...,10\}$. 
For each of 50 grid instances, two stages of trials are conducted, one which samples $\vec{C}$ from $[\vec{0},\vec{\infty}]$, keeping $\vec{D} = 0$, and the other which samples $\vec{D}$ from $[\vec{0},\vec{\infty}]$, keeping $\vec{C} = 0$. 
All components of the hyperparameter vectors are equal (this can be done because the cost functions are uniformly sampled from the same range). 
The reason for keeping one hyperparameter equal to zero while the other is modified is because we aim to tune between pure runtime efficiency (EMOA*) and pure memory efficiency (PIDMOA*), and because a high $\vec{D}$ reduces the portion of search performed using best-first search, the memory-efficiency gained by keeping $\vec{C} = \vec{0}$ outweighs the runtime benefit of a higher $\vec{C}$. This is especially true given the results we found below, which show that a high $\vec{D}$ sacrifices far more runtime compared to a low $\vec{C}$.
For each $k$ we tested, we chose a different set of $(\vec{C}, \vec{D})$ tuples because the choice of $\vec{C}$ and $\vec{D}$ depends on the branching factor and maximum depth of a solution, both of which vary between different $2^k$-connected grids. 
We found setting $\vec{C} > \max_e\vec{c}(e)$ to result in minimal memory improvement compared to EMOA*. This is because, for any child $l'$ of $l$, $\vec{0} \leq \vec{f}(l') - \vec{f}(l) = \vec{h}(v(l')) + \vec{c}(v(l),v(l')) - \vec{h}(v(l)) \leq \vec{c}(v(l),v(l'))$. Thus setting $\vec{C}$ any higher than $\max_e\vec{c}(e)$ did not generally result in any less children of $l$ being added to OPEN on first expansion of $l$ than EMOA*. In choosing $\vec{D}$, we tested a wide range of $\vec{D}$ from $\vec{0}$ to $\vec{\infty}$ and plotted those instances that were computationally feasible within a time constraint of 60 seconds. Because the optimal choice of $\vec{C}$ and $\vec{D}$ is problem dependent, there is no easy formulaic method of finding the best hyperparameters, and therefore trial and error worked most efficiently for us, especially for $\vec{D}$.

In Fig.~\ref{fig:grid}, we compare the performance of \alg{} against EMOA* for $2^2$- and $2^5$-connected grids. To note, while $\vec{D} = 0$, we see negligible improvements in memory usage over EMOA* when the branching factor $b$ is small ($k = 2 \implies b = 2^3 = 8$). However, when $k = 5$ ($b = 2^5 = 32$), we see a significant improvement in memory usage (down to 24.23\% of the memory usage of EMOA* on average) using \alg{} with $\vec{C} = \vec{0}$. Upon further increasing $\vec{D}$ to $\vec{16}$, \alg{} consumes only 5.03\% of the memory required of EMOA*. More notably, for $k = 2$, increasing $\vec{D}$ to $\vec{70}$ manages to reduce memory consumption to some degree (75.16\% of the memory consumption of EMOA*), where partial expansion alone ($\vec{C} = \vec{0}$ and $\vec{D} = \vec{0}$) could not. For both branching factors, we see non-negligible increases in runtime (88.36\% longer than EMOA* on average for $k = 5$) when using $\vec{C} = \vec{0}$ and $\vec{D} = \vec{0}$, and orders of magnitude worse runtime for large $\vec{D}$. However, for $\vec{C} = \vec{3}$ and $\vec{D} = \vec{0}$ (not plotted for visibility), we see substantially improved runtime over $\vec{C} = \vec{0}$ (only 8.06\% longer than EMOA* on average for $k = 5$) with nearly the same amount of memory usage (only 7.00\% increase in memory usage over $\vec{C} = \vec{0}$ on average for $k=5$). This implies that, for large branching factor problems, a well-chosen $\vec{C}$ can effectively reduce memory usage (to nearly the limit of the capabilities of \alg) without sacrificing significant runtime performance. It also shows that, if runtime is not as high a priority as memory efficiency, $\vec{D}$ can be used to further improve memory efficiency, at a heavy runtime cost.

\subsection{Robot State Lattice Test}\label{evaluation:lattice}
\begin{figure*}[t]
    \centering
    \begin{tabular}{c|c}
    \toprule
    \multicolumn{2}{c}{20x20 Robot State Lattice Test}\\
    \midrule
    $M = 2$ & $M = 3$\\
    \includegraphics[height=0.22\textwidth]{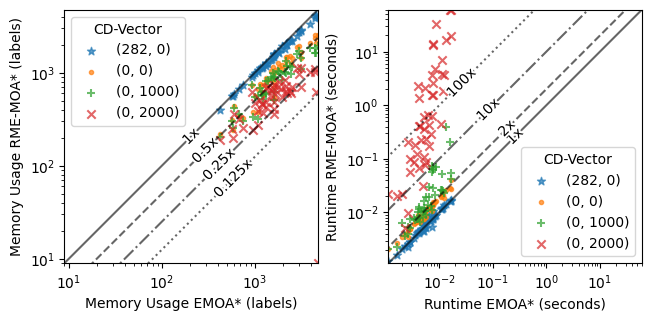} & \includegraphics[height=0.22\textwidth]{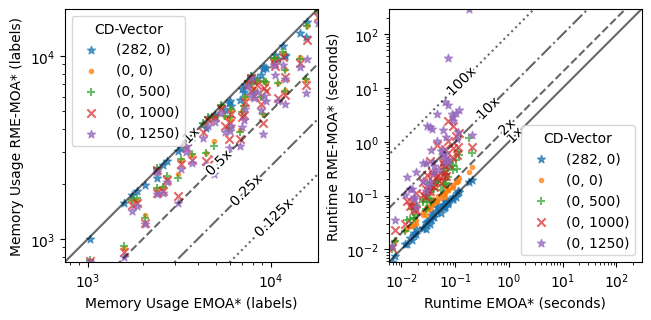}\\
    \end{tabular}
    \caption{Comparison of memory usage and runtime metrics for \alg{} with varying tuples $(\vec{C}, \vec{D})$, against EMOA* $(\vec{C}=\vec{\infty}, \vec{D}=\vec{0})$. The search space is a 20x20 robot state lattice with cost dimension $M=2$ (only distance and complexity objectives) and $M=3$ (adding a safety objective). These cost objectives are compared above ($M=2$ on the left, $M=3$ on the right). We create 50 randomized obstacle grid instances with obstacle density of 20\%, and we run \alg{} with all $(\vec{C}, \vec{D})$ tuples and EMOA*. Plot semantics are the same as in Fig.~\ref{fig:grid}.}
    \label{fig:lattice}
\end{figure*}
We then compare \alg{} against EMOA* on a mobile robot state lattice, which is similar to the $2^k$-connected grid test in Sec.~\ref{evaluation:grid_map}, with the exception that edges between vertices are ``differentially constrained''. This means that a robot at vertex $v_1$ can only move immediately to another $v_2$ if $v_2$ is reachable using at least one of a defined set of motion primitives (we refer the reader to~\cite{lattice} for a visualization of these primitives). Each vertex now represents a robot pose, discretized along $x$, $y$, and $theta$. Specifically, there exist 8 discrete angles a robot can take on at any ($x, y$) location. This results in a set of approximately 15 distinct motion primitives, corresponding to an average branching factor $b$ of 15. We also randomly introduce obstacles at $(x,y)$ locations with uniform density 0.2. We define 3 practically-oriented cost functions for search on this graph: path length, integrated turning angle over the course of the path, and safety of the path, which we define to be the number of obstacles within an 8-connected vicinity of a vertex.

In Fig.~\ref{fig:lattice}, we compare the performance of \alg{} against EMOA* on this robot state lattice. We see about half the memory consumption using $\vec{C} = \vec{0}$ and $\vec{D} = \vec{0}$ over EMOA* for $M=2$ (56.98\% of EMOA*), and nearly double the runtime (82.20\% increase). For $M=3$, we see even less of a memory improvement (78.99\% of EMOA*) but about the same runtime sacrifice (74.51\% increase). When $\vec{D}$ is raised until search becomes runtime infeasible, we see a small further improvement for $M=2$ and $M=3$ (30.67\% and 67.38\% of EMOA*, respectively), while runtime rises again by several orders of magnitude. We attribute this less impressive performance to a smaller branching factor in the robot state lattice ($\sim15$) than the grid map ($2^5=32$). We also believe that a lack of conflicting objectives may lead only a few promising partial paths to be continually expanded, resulting in a large $\alpha$ at each vertex in those paths, and thus a lack of efficacy of the partial expansion technique. This would also explain why \alg{} seems to perform worse on higher-dimensional cost functions.

\section{Conclusion}
This work develops a memory-efficient extension of the recent advanced runtime-efficient MOA* algorithms by generalizing the notion of partial expansion to the multi-objective case, integrating it with the runtime improvements of EMOA*, and fusing this algorithm with PIDMOA* to further improve memory efficiency. Our method, \alg, balances runtime and memory efficiency via two tunable hyperparameters, which can be set to implement the behavior of EMOA*, PIDMOA*, or some mixture of the two that balances runtime performance with memory efficiency. The results indicate that a properly tuned \alg{} is most effective at reducing memory usage for large branching factor problems with minimal runtime increases. \alg{} is less effective at reducing memory usage for small branching factor problems, where the ratio between OPEN and CLOSED is very small and partial expansion is less viable. For these problems, more memory efficiency can be gained if necessary by leaning more heavily on PIDMOA* search, at the expense of runtime performance.

For future work, the proposed approach in this work can be possibly further combined with other multi-objective search for both single-agent~\cite{BOBA*,Zhang2022Apex} and multi-agent~\cite{ren2022conflict} to balance between memory usage and runtime.

\section*{Acknowledgments}
This material is based upon work supported by the National Science Foundation under Grant No. 2120219 and 2120529. Any opinions, findings, and conclusions or recommendations expressed in this material are those of the author(s) and do not necessarily reflect the views of the National Science Foundation.
    
\bibliographystyle{IEEEtran}
\bibliography{rmemoa}

\end{document}